\documentclass[12pt]{amsart}
\usepackage{graphicx}
\usepackage[centertags]{amsmath}
\usepackage{amsfonts}
\usepackage{amssymb}
\usepackage{amsthm}
\usepackage{hyperref}
\newtheorem{thm}{Theorem}%[section]
\newtheorem{lem}{Lemma}[section]
\newtheorem{cor}{Corollary}
\newtheorem*{conj}{Conjecture}
\newtheorem{prop}[lem]{Proposition}
\theoremstyle{definition}

\theoremstyle{remark}
\newtheorem{rem}{Remark}[section]
\numberwithin{equation}{section}

\newcommand{\norm}[1]{\left\Vert#1\right\Vert}

\newcommand{\set}[1]{\left\{#1\right\}}

\newcommand{\T}[1]{\tilde{T}_{#1}}
\newcommand{\TT}[1]{T_{#1}}

\newcommand{\frakO}{\mathfrak{O}}

\newcommand{\calG}{\mathcal{G}}

\newcommand{\calN}{\mathcal{N}}

\newcommand{\calP}{\mathcal{P}}

\newcommand{\bbZ}{\mathbb{Z}}

\newcommand{\bbQ}{\mathbb{Q}}
\newcommand{\bbR}{\mathbb R}

\newcommand{\bbN}{\mathbb N}
\newcommand{\bbT}{\mathbb{T}}

\newcommand{\Tr}{ \mbox{Tr}}
\newcommand{\Op}{ \operatorname{Op}}
\newcommand{\SL}{ \mathrm{SL}}

\newcommand{\SU}{ \mathrm{SU}}

\newcommand{\Mat}{\mathrm{Mat}}

\newcommand{\Sq}{\mathrm{Sq}}

\newcommand{\Prob}{\mathrm{Prob}}
\newcommand{\calH}{\mathcal{H}}

\begin{document}
\title[cat map modulo prime powers]
{On matrix elements for the quantized\\ cat map modulo prime powers}%
\author{Dubi Kelmer }%
\address{School of Mathematics,
Institute for Advanced Study, 1 Einstein Drive , Princeton, New Jersey 08540 US}
\email{kelmerdu@ias.edu}

\thanks{
This material is based upon work supported by the National Science Foundation under agreement No. DMS-0635607. Any opinions, findings and conclusions or recommendations expressed in this material are those of the author and do not necessarily reflect the views of the National Science Foundation.}%

\date{\today}%
\dedicatory{}%
\commby{}%

\begin{abstract}
he quantum cat map is a model for a quantum system with underlying
chaotic dynamics. In this paper we study the matrix elements of
smooth observables in this model, when taking arithmetic symmetries
into account. We give explicit formulas for the matrix elements as
certain exponential sums. With these formulas we can show that there
are sequences of eigenfunctions for which the matrix elements decay
significantly slower then was previously expected. We also prove a
limiting distribution for the fluctuation of the normalized matrix
elements around their average.
\end{abstract}
\maketitle

\section{Introduction}
The quantum cat map is a model for a quantum system with underlying
chaotic dynamics that was originally introduced by the physicists
Hannay and Berry  \cite{HannayBerry80}. This model can be used to
study the semiclassical properties of such systems
\cite{BonechiDebievre,FaureNonnenmacherDeBievre03,KurRud2000,KurRud05}.
The classical dynamics underlying this model is the discrete time
iteration of a hyperbolic map, $A\in\SL(2,\bbZ)$, on the torus,
$\bbT^2=\bbR^2/\bbZ^2$. In order to quantize the cat map, for every
integer $N$ (playing the role of the inverse of Planck's constant)
the Hilbert space of states is $\calH_N=L^2(\bbZ/N\bbZ)$. For every
smooth real valued function $f$ there is a quantum observable, i.e.,
a Hermitian operator $\Op_N(f):\calH_N\to\calH_N$. The quantum
evolution is given by a unitary operator $U_N(A)$ on $\calH_N$.

For generic quantum systems with underlying chaotic dynamics, it is
believed that matrix elements of smooth observables tend to the
phase space average of the observable in the semiclassical limit. In
order to test this phenomenon in the quantum cat map  model,
Kurlberg and Rudnick introduced hidden symmetries of this model, a
group of commuting operators that commute with $U_N(A)$, they called
Hecke operators ~\cite{KurRud2000}. They showed that for any
sequence of Hecke eigenfunctions (i.e., joint eigenfunctions of all
Hecke operators), the corresponding matrix elements converge to the
phase space average as $N\to\infty$. To be more precise they showed
\cite[Theorem 1]{KurRud2000} that for any $f\in C^\infty(\bbT^2)$
and $\psi\in\calH_N$ a Hecke eigenfunction the matrix elements
satisfy
\[|\langle\Op_N(f)\psi,\psi\rangle-\int_{\bbT^2}f|\ll_{f,\epsilon} N^{-\frac{1}{4}+\epsilon}.\]
They remarked \cite[Remark 1.2]{KurRud2000} that the exponent of
$\frac{1}{4}$ is not optimal and that the correct bound should be
$O(N^{-\frac{1}{2}+\epsilon})$, in accordance to the second and
fourth moments. For $N$ prime (and consequently also for $N$ square
free) this is indeed the correct bound
\cite{DegliEspostiGraffiIsola95,GurevichHadani06}.

\begin{rem}
We note that without the arithmetic symmetries these bounds hold
only if the spectral degeneracies are sufficiently small. In fact,
there are sequences of eigenfunctions (where the degeneracies are
exceptionally large) that don't converge to the phase space average
at all. For these eigenfunctions the matrix elements localize around
short periodic orbits in the sense that the coresponding limiting
measure contains a component that is supported on the periodic orbit
\cite{FaureNonnenmacherDeBievre03}.
\end{rem}

In \cite{KurRud05} Kurlberg and Rudnick went on to investigate the
fluctuation of the normalized matrix elements,
\begin{equation}\label{e:F}
F_j^{(N)}=\sqrt{N}\left(\langle
\Op_N(f)\psi_j,\psi_j\rangle-\int_{\bbT^2}fdx\right),
\end{equation}
where $\psi_j$ are Hecke eigenfunctions and $N\to\infty$ through
primes. For this purpose they introduced the quadratic form
$Q(n)=\omega(nA,n)$ (with $\omega(n,m)=n_1m_2-n_2m_1$ the standard
symplectic form) and used it to define twisted Fourier coefficients.
For a smooth function $f\in C^\infty(\bbT^2)$ with Fourier
coefficients $\hat{f}(n)$ for $n\in\bbZ^2$, the twisted coefficients
are given by
\begin{equation}\label{e:Fourier}
f^{\#}(\nu)=\sum_{Q(n)=\nu}(-1)^{n_1n_2}\hat f(n).
\end{equation}

\begin{conj}[Kurlberg-Rudnick \cite{KurRud05}]
As $N\to\infty$ through primes, the limiting distribution of the
normalized matrix elements $F_j^{(N)}$ is that of the random
variable
\[X_f=\sum_{\nu\neq 0}
f^\#(\nu) \Tr(U_\nu)\] where $U_\nu$ are independently chosen random
matrices in $\SU(2)$ endowed with Haar probability measure.
\end{conj}
As evidence, the second and fourth moment were computed to show
agreement with this conjecture. In particular, the moment
calculation implies that the limiting distribution is not Gaussian,
in contrast to generic chaotic systems where the fluctuations are
believed to be Gaussian \cite{EckhardtFishman95,FeingoldPeres86}.

In this paper we further study the matrix elements for the cat map
for composite $N$. In fact, it is sufficient to understand the case
of prime powers (see \cite[Section 4.1]{KurRud2000}), and so we
restrict ourselves to this case. For $N$ a prime power, we give an
explicit formula for the matrix elements as a weighted sum of
certain exponential sums. We then use this formula to show that
there are sequences of eigenfunctions such that the matrix elements
decay like $N^{-1/3}$ rather then the expected rate of
$N^{-1/2+\epsilon}$. We further show that when $N=p^k$ with $k>1$,
the matrix elements have a limiting distribution as $p\to\infty$.
This distribution is not Gaussian and it is also different from the
(conjectured) distribution for $k=1$. Instead of behaving like
traces of random elements from $\SU(2)$, here the normalized matrix
elements vanish for half of the eigenfunctions and for the rest they
behave like $2\cos(\theta)$ where the angle is chosen at random.

\subsection{Results}
For every $N=p^k$ denote by
$$C(p^k)=\set{B\in \SL(2,\bbZ/p^k\bbZ)| AB=BA\pmod {p^k}},$$
the group of Hecke operators. For $\nu\in\bbZ$ and $\chi$ a
character of $C(p^k)$ define the exponential sum
\[E_{p^k}(\nu,\chi)=\sum_{x\in X(p^k)}e_{p^k}(\nu x)\chi(\beta(x)),\]
where $$X(p^k)=\set{x\in\bbZ/p^k\bbZ|(\Tr(A)^2-4)x^2\neq
1\pmod{p}}$$ and $\beta\colon X(p^k)\hookrightarrow C(p^k)$ is an
injection of $X(p^k)$ into $C(p^k)$ given by a rational function
(defined by (\ref{e:beta})).

\begin{thm}\label{t:form1}
For each prime power $p^k$, there is a subset $\hat C_0(p^k)\subset
\hat{C}(p^k)$ of characters, with $\lim_{p\to\infty}\frac{|\hat
C_0(p^k)|}{p^k}=1$ such that
\begin{enumerate}
\item For any $\chi\in \hat C_0(p^k)$ there is a unique Hecke eigenfunction $\psi$, s.t., $\chi$ is a joint eigenvalue.
\item For this eigenfunction, and any elementary observable $f_n(x)=\exp(2\pi i n\cdot x)$ with $Q(n)\not\equiv 0\pmod{p}$
\[\langle\Op_{p^k}(f_n)\psi,\psi\rangle=\pm\frac{(-1)^{n_1n_2}}{\# C(p^k)}E_{p^k}(\frac{Q(n)}{2},\chi\chi_0),\]
where $\chi_0$ is a fixed character of $C(p^k)$ and the sign $(\pm)$
depends on $p,k$ but not on $\psi$.
\end{enumerate}
\end{thm}

If we consider nontrivial prime powers (i.e., $k>1$) we can use
elementary methods to evaluate these sums. In particular we find
that there are matrix elements that decay much slower then the
expected rate of $N^{-\frac{1}{2}+\epsilon}$.
\begin{thm}\label{t:slowdecay}
There are smooth observables $f\in C^\infty(\bbT^2)$, and sequences
of Hecke eigenfunctions satisfying
$|\langle\Op_N(f)\psi_j,\psi_j\rangle-\int_{\bbT^2} f|\gg
N^{-\frac{1}{3}}$.
\end{thm}

We note, however, that these exceptional matrix elements are quite
rare, in the sense that for a fixed observable the number of matrix
elements decaying slower then $N^{-\frac{1}{2}+\epsilon}$ is bounded
by $O(p^{k-1})$ (see Corollary \ref{c:bound}) .

\begin{rem}
In \cite{Olofsson07}  Olofsson studied the supremum norm of Hecke
eigenfunctions for the quantized cat map. He showed that for
composite $N$ the supremum norm can be of order $N^{\frac{1}{4}}$,
which is much larger then the case of $N$ prime (or square free)
where all Hecke eigenfunctions satisfy $\norm{\psi}_\infty\ll
N^\epsilon$ \cite{GurevichHadanisup05,Kurlberg07}. Although the two
phenomena look similar, there does not seem to be any apparent
connection between them. At least in the sense that the
eigenfunctions with large matrix elements are usually not the
eigenfunctions with large supremum norm.
\end{rem}

For nontrivial prime powers, we can also show that the exponential
sums $E_{p^k}(\nu,\chi)$ (and hence also the matrix elements) have a
limiting distribution as $p\to\infty$. (See \cite{KelmerSums08} for
similar results on twisted Kloosterman sums). To simplify the
discussion we will assume from here on that the observable $f$ is a
trigonometric polynomial and let $F_j^{(N)}$ be the normalized
matrix element as in (\ref{e:F}). Let $\mu$ denote the measure on
$[0,\pi)$ defined by
\begin{equation*}
\mu(f)=\frac{1}{2}f(\frac{\pi}{2})+\frac{1}{2\pi}\int_0^{\pi}f(\theta)d\theta.
\end{equation*}

\begin{thm}\label{t:main}
Let $f$ be a trigonometric polynomial. For any $k>1$, as $p\to
\infty$ through primes, the limiting distribution of the normalized
matrix elements $F_j^{(p^k)}$ is that of the random variable
\[Y_f=2\sum_{\nu\neq 0}f^\#(\nu)\cos(\theta_\nu)\]
where $\theta_\nu$ are independently chosen from $[0,\pi)$ with
respect to the measure $\mu$.
\end{thm}
\begin{rem}
As mentioned above, there can be exceptionally large matrix elements
for which $F_j^{(N)}\gg N^{1/6}$ are not bounded. Such matrix
elements would cause the moments (above the $6$'th moment) to blow
up as $N\to\infty$. Nevertheless, since the number of exceptional
matrix elements is of limiting density zero, they do not influence
the limiting distribution (see section \ref{ss:equi} for more
details).
\end{rem}
\subsection{Outline}
The outline of the paper is as follows: In section \ref{s:BKG} we
provide some background on the cat map and its quantization and on
the notion of a limit distributions. In section \ref{s:form} we
compute the formulas for the matrix elements proving Theorem
\ref{t:form1}. In section \ref{s:comp} we compute the exponential
sums appearing in these formulas for non trivial prime powers, and
establish the limiting distribution as $p\to\infty$. Then in section
\ref{s:equi} we deduce both of the results on the matrix elements
(Theorems \ref{t:slowdecay} and \ref{t:main}) from the analysis of
the exponential sums.
\section{Background}\label{s:BKG}
The full details for the cat map and it's quantization can be found
in ~\cite{KurRud2000}. We briefly review the setup and go over our
notation.

\subsection{Classical dynamics}
The classical dynamics are given by the iteration of a hyperbolic
linear map $A\in \SL(2,\bbZ)$.
\[x=\left(%
\begin{array}{c}
p \\
q \\
\end{array}%
\right)\in\mathbb{T}^2\mapsto Ax\pmod{1}.\] Given an observable
$f\in C^\infty(\mathbb{T}^2)$, the classical evolution defined by
$A$ is $f\mapsto f\circ A$.
\subsection{Quantum kinematics}
For doing quantum mechanics on the torus, one takes Planck's
constant to be $1/N$, as the Hilbert space of states one takes
$\mathcal{H}_N=L^2(\bbZ/N\bbZ)$, where the inner product is given
by:
\[\langle\phi,\psi\rangle=\frac{1}{N}\sum_{y \in\bbZ/N\bbZ}\phi(y)\overline{\psi(y)}.\]
For $n=(n_1,n_2)\in \bbZ^2$ define elementary operators $\TT{N}(n)$
acting on $\psi\in \mathcal{H}_N$ via:
\begin{equation}\label{eTN:1}
\TT{N}(n)\psi(y)=e_{2N}(n_1n_2)e_N(n_2y)\psi(y+n_1),
\end{equation}
where $e_N(x)=e^{\frac{2\pi i x}{N}}$. For any smooth classical
observable $f\in C^\infty(\mathbb{T}^2)$ with Fourier expansion
$f(x)=\sum_{n\in\bbZ^2}\hat{f}(n)e^{2 \pi i n\cdot x}$, its
quantization is given by
\[\Op_N(f)=\sum_{n\in\bbZ^2}\hat{f}(n)T_N(n).\]

\subsection{Quantum dynamics:} For any $A\in \SL(2,\bbZ)$, we assign
unitary operators $U_N(A)$, acting on $L^2(\bbZ/N\bbZ)$ having the
following important properties:
\begin{itemize}
\item ``Exact Egorov'': For $A\equiv I\pmod{2}$, and any $f\in
C^\infty(\mathrm{\bbT}^2)$
\[U_N(A)^{-1}\Op_N(f)U_N(A)=\Op_N(f\circ A).\]
\item The map $A\mapsto U_N(A)$ is a representation of $\SL(2,\bbZ/N\bbZ)$:
If $C\equiv AB\pmod{N}$ then $U_N(A)U_N(B)=U_N(C)$.
\end{itemize}

We will make use of the following formula for $U_N(A)$, (valid for
odd $N$ and any $A\in\SL(2,\bbZ)$) \cite[Proposition
1.4]{Kelmer05Arith}.
\begin{equation}\label{e:form}
U_{N}(A)=\frac{\sigma_N(A)}{|\ker_N(A-I)|N}\sum_{m\in (\bbZ/N
\bbZ)^2}\T{N}(m)\T{N}(-mA)
\end{equation}
where $\sigma_N(A)=\Tr(U_N(A))$ is the character of the
representation,
$$|\ker_N(A-I)|=\#\set{n\in(\bbZ/N\bbZ)^2|n(A-I)\equiv 0\pmod{N}},$$
and $\T{N}(n)=(-1)^{n_1n_2}\TT{N}(n)$ are  twisted elementary
operators.
\begin{rem}
The twisted operators $\T{N}(n)$ have the convenient feature that
$$U_N(A)^*\T{N}(n)U_N(A)=\T{N}(nA)$$ for any $A\in\SL(2,\bbZ/N\bbZ)$ (without the parity condition).
\end{rem}

\subsection{Hecke eigenfunctions}
Let $\alpha>\alpha^{-1}$ be the eigenvalues of $A$ in a (real)
quadratic extension $K/\bbQ$. Then the vectors
$\vec{v}_\pm=(c,\alpha^{\pm 1}-a)$ are corresponding eigenvectors
$\vec{v}_\pm A=\alpha^{\pm 1}\vec{v}_{\pm}$. Denote by
$D=\Tr(A)^2-4\in\bbZ^{+}$ so that $\sqrt{D}=\alpha-\alpha^{-1}$.
Consider the ring $\frakO=\bbZ[\alpha]$ and denote by
$\iota:\frakO\to\Mat(2,\bbZ)$ the map sending
$\beta=n+m\alpha\mapsto B=n+mA$ (this map is a ring homomorphism as
$\alpha$ and $A$ have the same minimal polynomial).

For any integer $N$ the norm map, $\calN_{K/\bbQ}:K^*\to\bbQ^*$,
induce a well defined map
$\calN_{N}:(\frakO/N\frakO)^*\to(\bbZ/N\bbZ)^*$. Let $C(N)=\ker
\calN_N$ be its kernel, then its image $\iota(C(N))\subset
\SL(2,\bbZ/N\bbZ)$ is a commutative subgroup of $\SL(2,\bbZ/N\bbZ)$
that commutes with $A \pmod{N}$. The Hecke operators are then
$\set{U_N(B)|B\in\iota(C(N))}$, and Hecke eigenfunctions are joint
eigenfunctions of $U_N(A)$ and all the Hecke operators.

The eigenvalues corresponding to each Hecke eigenfunction define a
character $\chi$ of $C(N)$ i.e.,
$U_N(\iota(\beta))\psi=\chi(\beta)\psi$. We can thus decompose our
Hilbert space into a direct sum of joint eigenspaces
$\calH_N=\bigoplus_{\chi} \calH_\chi$, parameterized by the
characters of $C(N)$. We say that a character $\chi$ appears with
multiplicity one in the decomposition when the corresponding
eigenspace is one dimensional.

\subsection{Limit distribution}\label{ss:equi}
    We recall the notion of a limiting distribution for a sequence of points on the line.
    For each $N$ let $\{F_j^{(N)}\}_{j=1}^N$ be a set of points on the line. We say
    that these points have a limiting distribution $Y$ (where $Y:\Omega\to\bbR $ is some
    random variable on a probability space $\Omega$) if for any segment $[a,b]\subset\bbR$ the limit
    \[\lim_{N\to\infty}\frac{\#\set{j|a\leq F_j^{(N)}\leq
    b}}{N}=\Prob(Y\in[a,b]).\]
    From this definition it is strait forward that making an arbitrary change
    in a density zero set of points
    (i.e., changing $S_N$ points for each $N$ with
    $\frac{S_N}{N}\to 0$),
    does not affect the limiting distribution.

    An equivalent condition for having a limiting distribution $Y$,
    is that for any continues bounded
    function $g$ the average $\frac{1}{N}\sum_j g(F_j^{(N)})$
    converges as $N\to\infty$ to $\int_\Omega g(Y(\omega))d\omega$.
    Note that the condition that the test function $g$ is
    bounded is necessary unless both the variable $Y$ and the points $F_j^{(N)}$
    are uniformly bounded. In particular, if the points $F_j^{(N)}$
    are not uniformly bounded then their moments don't necessarily converge to the
    moments of $Y$.

\subsection{Notation}
We use the notation $e(x)=e^{2\pi i x}$. For any $N\in\bbN$ we
denote by $e_N(\cdot)$ the character of $\bbZ/N\bbZ$ given by
$e_N(x)=e(\frac{x}{N})$. When there is no risk of confusion we will
slightly abuse notation and write $e_N(\frac{a}{b})$ for
$e_N(ab^{-1})$ (where $b^{-1}$ denotes the inverse of $b$ modulo
$N$). For example, for $N$ odd and $a\in\bbZ$ we may write
$e_{2N}(a)=(-1)^ae_N(\frac{a}{2})$.

\section{Formulas for Matrix Elements}\label{s:form}
For $N$ a prime power we give formulas for the matrix elements of
elementary observables explicitly as exponential sums. When $N$ is
prime these formulas appeared in \cite{KurRud05} (for primes that
split in $\frakO$) and in \cite{Kelmer05Arith} (for inert primes).

We will make use of the following parametrization of the Hecke
operators. For any integer $1\leq l\leq k$ we define subgroups
$C_p(k,l)\subset C(p^k)$ by
\[C_p(k,l)=\set{\beta\in C(p^k)|\beta\equiv 1\pmod {p^l}}.\]
For notational convenience we will also define $C_p(k,k+1)=\{1\}$.
Let $$X(p^k)=\set{x\in\bbZ/p^k\bbZ|Dx^2\neq 1\pmod{p}}$$ then the
map
\begin{equation}\label{e:beta}
\beta(x)=\frac{\sqrt{D}x+1}{\sqrt{D}x-1}
\end{equation}
 is a bijection between $X(p^k)$ and $C(p^k)\setminus C_p(k,1)$ with inverse map given by $x=\frac{1+\beta(x)}{\sqrt{D}(1-\beta(x))}\pmod{p^k}$ (note that for $\beta\neq 1\pmod{p}$ the inverse map is indeed well defined).
For every character $\chi$ of $C(p^k)$ and any $\nu\in
(\bbZ/p^k\bbZ)^*$ we have the exponential sum
\[E_{p^k}(\nu,\chi)=\sum_{x\in X(p^k)}e_{p^k}(\nu x)\chi(\beta(x)).\]
To prove Theorem \ref{t:form1} we will show that for any
$n\in\bbZ^2$ with $Q(n)=\nu\neq 0\pmod{p}$, and for every character
$\chi$ of $C(p^k)$ that appears with multiplicity one, the
corresponding matrix element is given by
 \[\langle \T{p^k}(n)\psi,\psi\rangle=\frac{\pm 1}{\# C(p^k)}\sum_{x\in X(p^k)}e_{p^k}(\frac{\nu x}{2})\chi\chi_0(\beta(x)),\]
(where $\chi_0$ is a fixed character of $C(p^k)$ and the sign is
$-1$ when $p$ is inert and $k$ is odd and $+1$ otherwise). We can
then take our set $\hat{C}_0(p^k)$ to be the set of characters
appearing with multiplicity one. This set is of order $p^k$ if $p$
is inert (Lemma \ref{l:mult}) and of order $p^{k}-p^{k-1}$ if $p$
splits (Lemma \ref{l:eigenfunction}). Hence,
 indeed $\frac{\hat{C}_0(p^k)}{p^k}=1+O(\frac{1}{p})$.
We will compute the matrix elements separately for the inert and
split cases.

\subsection{Split case}
When $p$ is split, we can give explicit formulas for the Hecke
eigenfunctions and use them to compute the matrix elements. Since we
assume that $p$ splits in $\frakO$, there is a matrix
$M\in\SL(2,\bbZ/p^k\bbZ)$ satisfying that $M^{-1}AM=\begin{pmatrix}
y & 0\\ 0 & y^{-1}\end{pmatrix}\pmod{p^k}$. Consequently, the Hecke
group is given by
$$C(p^k)=\set{M\begin{pmatrix}x & 0\\ 0 & x^{-1}\end{pmatrix}M^{-1}| x\in(\bbZ/p^k\bbZ)^*},$$ which is naturally isomorphic to $(\bbZ/p^k\bbZ)^*$. We recall that
\begin{equation}\label{eq:diag}
U_{p^k}(\begin{pmatrix} x & 0\\ 0 &
x^{-1}\end{pmatrix})\psi(y)=\chi_0(x)\psi(xy),
\end{equation}
where $\chi_0$ is a fixed character of $(\bbZ/p^k\bbZ)^*$
\cite[Section 4.3]{KurRud2000}.

\begin{lem}\label{l:eigenfunction}
For any character $\chi$ of $(\bbZ/p^k\bbZ)^*$ (extended to a
function on $\bbZ/p^k\bbZ$ by setting $\chi(px)=0$), the function
$\psi=\sqrt{\frac{p}{p-1}}U_{p^k}(M)\chi$ is a normalized joint
eigenfunction of all Hecke eigenfunctions with eigenvalue
$\chi\chi_0$. Furthermore, if $\chi$ is not trivial on the subgroup
$C_p(k,k-1)$ then this is the only eigenfunction.
\end{lem}
\begin{proof}
The first assertion is an immediate consequence of (\ref{eq:diag}).
For the second part, assume that $\psi$ is an eigenfunction with
eigenvalue $\chi\chi_0$, and that $\chi$ is not trivial on
$C_p(k,k-1)$. Then there is $x_0\in C_p(k,k-1)$ with $\chi(x_0)\neq
1$. Now, let $\phi=U_{p^k}(M)^{-1}\psi$, then for any
$x\in(\bbZ/p^k\bbZ)^*$,
\[\chi\chi_0(x)\phi(y)=U_{p^k}(\begin{pmatrix}x & 0\\ 0 & x^{-1}\end{pmatrix})\phi(y)=\chi_0(x)\phi(xy),\]
hence $\phi(xy)=\chi(x)\phi(y)$. For any $y\equiv 0\pmod{p}$ we have
that $x_0y\equiv y\pmod{p^k}$ (as $x_0\equiv 1\pmod {p^{k-1}}$).
Consequently, $\phi(y)=\phi(x_0y)=\chi(x_0)\phi(y)$ implying that
$\phi(y)=0$. On the other hand, for $y\not\equiv 0\pmod{p}$ we have
$\phi(y)=\chi(y)\phi(1)$ so $\phi$ is uniquely determined (up to
normalization).
\end{proof}

\begin{rem}
In the case that the character $\chi$ is trivial on the group
$C_p(k,l)$ (but not on $C_p(k,l-1)$) then the above argument implies
that the corresponding eigenspace is of dimension $k-l+1$.
\end{rem}

\begin{proof}[\bf{Proof of Theorem \ref{t:form1}} (split case)]
Let $\chi$ be a character not trivial on $C_p(k,k-1)$. Then,
$\psi=\sqrt{\frac{p}{p-1}}U_{p^k}(M)\chi$ is an eigenfunction with
character $\chi\chi_0$, where $A=M\begin{pmatrix} y & 0\\ 0 &
y^{-1}\end{pmatrix}M^{-1}\pmod{p^k}$. Consequently, for any
(twisted) elementary observable
\begin{eqnarray*}
 \langle\T{p^k}(n)\psi,\psi\rangle &=&\frac{p}{p-1}\langle U_{p^k}(M)^*\T{p^k}(n)U_{p^k}(M)\chi,\chi\rangle\\
&=&\frac{p}{p-1}\langle \T{p^k}(m)\chi,\chi\rangle
\end{eqnarray*}
 with $m=nM\pmod{p^k}$.

Now, let $d=y-y^{-1}$ so that $d^2\equiv D \pmod {p^k}$ (recall
$\Tr(A)\equiv y+y^{-1}\pmod{p^k})$. Then
\begin{eqnarray*}
\langle \T{p^k}(m)\chi,\chi\rangle & = & \frac{1}{p^k}\sum_{x\in(\bbZ/p^k\bbZ)^*}e_{p^k}(\frac{m_1m_2}{2})e_{p^k}(m_2x)\chi(\frac{x+m_1}{x})\\
& = & \frac{1}{p^k}\sum_{t\in X(p^k)}e_{p^k}(\frac{dm_1m_2}{2}t)\chi(\frac{dt+1}{dt-1})\\
\end{eqnarray*}
where we made the change of variables $2x=m_1(dt-1)$. Finally,
notice that for $m=nM\pmod{p^k}$ we have that
\[Q(n)=\omega(nA,n)\equiv m_1m_2(y-y^{-1})\equiv dm_1m_2\pmod{p^k}.\]
Hence indeed
\begin{eqnarray*}
\langle \T{p^k}(n)\psi,\psi\rangle = \frac{1}{\# C(p^k)}\sum_{x\in
X(p^k)}e_{p^k}(\frac{Q(n)x}{2})\chi(\beta(x))
\end{eqnarray*}
\end{proof}

\subsection{Inert case}
First we show that for $p$ inert, any joint eigenspace is one
dimensional.

\begin{lem}\label{l:mult}
For $N=p^k$ and $p$ inert, the dimension of any joint eigenspace
satisfies $\dim \calH_\chi\leq 1$.
\end{lem}
\begin{proof}
The trace of the quantum propagators satisfy \cite[Corollary
1.6]{Kelmer05Arith}
\[|\Tr(U_{p^k}(B))|^2=\#\set{n\in(\bbZ/p^k\bbZ)^2:n(B-I)\equiv
0\pmod{p^k}}.\] For $p$ inert, the group $C(p^k)$ is of order $\#
C(p^k)=p^{k-1}(p+1)$, and the groups
\[C_p(k,l)=\set{\beta\in C(p^k)|\beta\equiv 1\pmod{p^l}},\]
 are of order $\# C_p(k,l)=\frac{\# C(p^k)}{\# C(p^l)}=p^{k-l}$.
Moreover, for any $\beta\in C_{p}(k,l)\setminus C_{p}(k,l+1)$ we
have
 $|\Tr(U_{p^k}(\iota(\beta)))|^2=p^{2l}$.
Consequently
 \begin{eqnarray*}
 \sum_{\beta\in C(p^k)}
 |\Tr(U_{p^k}(\iota(\beta)))|^2 & = & p^{k}+\sum_{l=1}^{k-1} \sum_{\beta\in C_p(k,l)\setminus
 C_p(k,l+1)}\!\!\!\!\!\!\!\!\!|\Tr(U_{p^k}(\iota(\beta)))|^2+p^{2k}\\
 & = & p^{k}+\sum_{l=1}^{k-1} (p^{k-l}-p^{k-l-1})p^{2l}+ p^{2k}\\
% & = & p^{k}+\sum_{l=1}^{k-1} (p^{k+l}-p^{k+l-1})+p^{2k}\\
 & = & p^{k}+p^{2k-1}-p^{k}+p^{2k}=p^k\#C(p^k)\\
\end{eqnarray*}
On the other hand, if we denote by $n_\chi=\dim\calH_\chi$ then
\[\frac{1}{\# C(p^k)}\sum_{\beta\in C(p^k)}
 |\Tr(U_{p^k}(\iota(\beta)))|^2=\sum_\chi n_\chi^2.\]
 Comparing the two expressions we get
 \[\sum_\chi n_\chi^2=\frac{1}{\# C(p^k)}\sum_{\beta\in C(p^k)}
 |\Tr(U_{p^k}(\iota\beta)))|^2=p^k=\dim\calH=\sum_\chi n_\chi.\]
 Since $n_\chi$ are non negative integers this implies $n_\chi\leq 1$.
\end{proof}

After establishing this fact, following the idea of Gurevich and
Hadani \cite{GurevichHadani06}, we can write the matrix elements of
elementary observables as
\[\langle \TT{p^k}(n)\psi_j,\psi_j\rangle=\Tr(\TT{p^k}(n)\calP_{\chi_j}),\]
with
\[\calP_{\chi_j}=\frac{1}{\# C(p^k)}\sum_{\beta\in C(p^k)}U_{p^k}(\iota(\beta))\bar\chi_j\]
the projection operator to the (one dimensional) eigenspace spanned
by $\psi_j$. We then use formula (\ref{e:form}) for
$U_{p^k}(\iota(\beta))$ in order compute
$\Tr(\TT{p^k}(n)U_{p^k}(\iota(\beta)))$. However, in order to do
this we first need to give a formula for the character of the
representation $\sigma(B)=\Tr(U_{p^k}(B))$ (appearing in
(\ref{e:form} )), for any $B\in\iota(C(p^k))$.

\begin{prop}
There is a character $\chi_0\in \hat C(p^k)$ such that for any
$\beta\in C(p^k)$, we have
\[\Tr(U_{p^k}(\iota(\beta)))=(-1)^{k}(-p)^l\chi_0(\beta),\]
with $1\leq l\leq k$ the maximal integer such that $\beta\equiv
1\pmod {p^l}$.
\end{prop}

\begin{proof}
For any $1\leq l\leq k+1$ consider the subgroup of characters
\[\hat C^{(l)}(p^k)=\set{\chi\in \hat C(p^k)|\chi(\beta)=1,\;\forall \beta\in C_p(k,l)}.\]
For $1\leq l\leq k$, the group $\hat C^{(l)}(p^k)$ is the kernel of
the restriction map from $\hat C(p^k)$ to $\hat C_p(k,l)$ and hence
of order $\frac{\# C(p^k)}{\# C_p(k,l)}=p^{l-1}(p+1)$ (and for
$l=k+1$ we have $\hat C^{(k+1)}(p^k)=\hat C^{(k)}(p^k)=C(p^k)$).

We will first prove the following: For each $1\leq l\leq k+1$ there
is a character $\chi_l\in \hat C^{(l)}(p^k)$ and a subset
$S_l\subset \hat C^{(l-1)}(p^k)$ of order $\# S_l=p^{l-1}$ such that
for any $\beta\in C(p^k)\setminus C_p(k,l)$,
\[\Tr(U_{p^k}(\iota(\beta)))=(-1)^{k+l+1}\chi_l\chi_{l+1}\cdots\chi_k(\beta)\sum_{\chi\in
S_l}\chi(\beta).\]

First for $l=k+1$ we take the character to be the trivial character
and the set $S_{k+1}\subset \hat C^{(k)}(p^k)=\hat C(p^k)$ to be the
set of characters that appear in the decomposition of $\calH_{p^k}$
(there are $p^k$ such characters each appearing with multiplicity
one). Then indeed $\Tr(U_{p^k}(\iota(\beta)))=\sum_{\chi\in
S_{k+1}}\chi(\beta).$ If $k=1$ the sum is over all but one of the
characters, say $\chi_0\in \hat C(p)$, and hence
$\Tr(U_{p}(\iota(\beta)))=-\chi_0(\beta)$ as claimed. For $k>1$ we
proceed by induction as follows.

We assume the assertion is true for $1<l\leq k+1$ and show that it
is true for $l-1$. For any $\beta\in C(p^k)\setminus
 C_p(k,l)$, by our assumption
\[\Tr(U_{p^k}(\iota(\beta)))=(-1)^{k+l+1}\chi_{l}\chi_{l+1}\cdots\chi_k(\beta)\sum_{\chi\in
S_{l}}\chi(\beta)\] with $S_l\subset \hat C^{(l-1)}(p^k)$ of order
$\# S_l=p^{l-1}$. The order $\# \hat C^{(l-1)}(p^k)=p^{l-1}+p^{l-2}$
hence the complement $S_l^c$ in $\hat C^{(l-1)}(p^k)$ is of order
$p^{l-2}$. Now, if $\beta\not\in C_p(k,l-1)$ then the sum over all
characters in $\hat C^{(l-1)}(p^k)$ vanish, and hence $\sum_{\chi\in
S_{l}}\chi(\beta)=-\sum_{\chi\in S_{l}^c}\chi(\beta)$. We thus have
that
\[\Tr(U_{p^k}(\iota(\beta)))=(-1)^{k+l}\chi_{l}\chi_{l+1}\cdots\chi_k(\beta)\sum_{\chi\in S_{l}^c}\chi(\beta).\]
On the other hand, for $\beta\in C_p(k,l-2)\setminus C_p(k,l-1)$ we
have that $|\Tr(U_{p^k}(\iota(\beta)))|=p^{l-2}$, which could happen
only if $\chi(\beta)$ takes the same value for all $\chi\in
S_{l}^c$. Now take $\chi_{l-1}$ to be any character from $S_{l}^c$
and let $S_{l-1}=\chi_{l-1}^{-1}S_{l}^c$. Then $S_{l-1}\subseteq
\hat C^{(l-2)}(p^k)$ is of order $p^{l-2}$ and
\[\Tr(U_{p^k}(\iota(\beta)))=(-1)^{k+l-1}\chi_{l-1}\chi_l\chi_{l+1}\cdots\chi_k(\beta)\sum_{\chi\in S_{l-1}}\chi(\beta).\]

Now, let $\chi_0=\chi_1\cdot\chi_2\cdots \chi_k$ and let $\beta\in
C_p(k,l)\setminus C_p(k,l+1)$. Since $\beta\not\in C_p(k,l+1)$ we
have,
\[\Tr(U_{p^k}(\iota(\beta)))=(-1)^{k+l}\chi_{l+1}\chi_{l+2}\cdots\chi_k(\beta)\sum_{\chi\in
S_{l+1}}\chi(\beta).\] On the other hand we also assume $\beta\in
C_p(k,l)$, hence, for all $\chi\in  S_{l+1}\subset \hat
C^{(l)}(p^k)$ we have $\chi(\beta)=1$ implying that $\sum_{\chi\in
S_{l+1}}\chi(\beta)=\# S_{l+1}=p^l$. Also for any $m\leq l$,
$\chi_m\in C^{(m)}(p^k)\subset C^{(l)}(p^k)$, so $\chi_m(\beta)=1$.
We thus get that indeed
\[\Tr(U_{p^k}(\iota(\beta)))=(-1)^k(-p)^l\chi_0(\beta).\]
\end{proof}

\begin{prop}\label{p:character}
Let $n\in\bbZ^2$ and $B \in \iota(C(p^k))$. For $B\equiv I\pmod{p}$
the trace $\Tr(\T{p^k}(n)U_{p^k}(B))=0$. Otherwise, there is $x\in
X(p^k)$ such that $B=\iota(\beta(x))$ and
\[\Tr(\T{p^k}(n)U_{p^k}(B))=(-1)^k
 \chi_0(\beta(x))e_{p^k}(-\frac{Q(n)x}{2}).\]
\end{prop}
\begin{proof}
Use formula (\ref{e:form}) for $U_{p^k}(B)$ to get that
\begin{eqnarray*}\lefteqn{\Tr(\T{p^k}(n)U_{p^k}(B))=}\\
&&\frac{\sigma_{p^k}(B)}{|\ker_{p^k}(B-I)|p^k}\sum_{m\in (\bbZ/p^k
\bbZ)^2}\Tr(\T{p^k}(n)\T{p^k}(m)\T{p^k}(-mB))
\end{eqnarray*}
Note that up to a phase
$\T{p^k}(n)\T{p^k}(m)\T{p^k}(-mB)=e^{i\alpha}\T{p^k}(n-m(B-I))$ and
recall that $\Tr(\T{p^k}(n))=0$ unless $n\equiv 0\pmod{p^k}$ (see
e.g., \cite[Lemma 4]{KurRud2000}). Hence, the only summand that does
not vanish is the one satisfying $n=m(B-I)\pmod {p^k}$. We can
assume $n\neq 0\pmod p$, so that the trace vanishes whenever
$B\equiv I\pmod{p}$. Otherwise, $B=\iota(\beta)$ for some $\beta\in
C(p^k)\setminus C_p(k,1)$ and $\sigma_{p^k}(B)=(-1)^k\chi_0(\beta)$
so that
\[\Tr(\T{p^k}(n)U_{p^k}(B))=(-1)^k\chi_0(\beta)e_{p^k}(-\frac{\omega(m,mB)}{2}),\]
with $m=n(B-I)^{-1}\pmod {p^k}$.

Now recall the parametrization $C(p^k)\setminus
C_p(k,1)=\set{\beta(x)|x\in X(p^k)}$, with
$\beta(x)=\frac{\sqrt{D}x+1}{\sqrt{D}x-1}$. We claim that for
$B=\iota(\beta(x))$ and $m=nB$ we have that $\omega(m,mB)=Q(n)x$. To
show this substitute $(\beta(x)-1)^{-1}=\frac{\sqrt{D}x-1}{2}$ and
$(\beta(x)-1)^{-1}\beta(x)=\frac{\sqrt{D}x+1}{2}$. Consequently we
get
\[\omega(m,mB)=\omega\left(n\iota(\frac{\sqrt{D}x-1}{2}),n\iota(\frac{\sqrt{D}x+1}{2})\right)=
\frac{x}{2}\omega(n\iota(\sqrt{D}),n).\] Recall that
$\sqrt{D}=(\alpha-\alpha^{-1})$ so that indeed
\[\omega(n\iota(\sqrt{D}),n)=\omega(n(A-A^{-1}),n)=2\omega(nA,n)=2Q(n).\]
\end{proof}

\begin{proof}[ \bf{Proof of Theorem \ref{t:form1}} (inert case)]
For every character $\chi$ let
\[\calP_\chi=\frac{1}{\# C(p^k)}\sum_{B\in C(p^k)} U_{p^k}(B)\bar\chi(B),\]
be the projection operator to the (one dimensional) eigenspace
corresponding to $\chi$. Let $\psi$ be the corresponding Hecke
eigenfunction. Then
\[\langle \T{p^k}(n)\psi,\psi\rangle=\Tr(\T{p^k}(n)\calP_\chi)=\frac{1}{\# C(p^k)}\sum_{B\in C(p^k)} \Tr(\T{p^k}(n)U_{p^k}(B))\bar\chi(B).\]
Now from the Proposition \ref{p:character}
\[\Tr(\T{p^k}(n)U_{p^k}(B))=(-1)^k\chi_0(B)e_{p^k}(-\frac{Q(n)x}{2}),\]
implying that
\[\langle \T{p^k}(n)\psi,\psi\rangle=\frac{(-1)^k}{\# C(p^k)}\sum_{\beta\in C(p^k)} e_{p^k}(-\frac{Q(n)x}{2})\chi_0(\beta(x))\bar\chi(\beta(x)).\]
After a change of variables $x\mapsto -x$ we get
\[\langle \T{p^k}(n)\psi,\psi\rangle=\frac{(-1)^k}{\# C(p^k)}E_{p^k}(Q(n)/2,\chi\bar\chi_0)\]
\end{proof}

\section{Analysis of the Exponential Sums}\label{s:comp}
In this section we compute the exponential sums $E_{p^k}(\nu,\chi)$
for any prime power $k>1$. This can be done using elementary methods
(see, e.g., \cite[section 12.3]{IwaniecKowalski} or \cite[Chapter
1.6]{BerndtEvansWilliams}), however, since the setup here is
slightly different we will perform this computation in full. We then
evaluate all mixed moments of these exponential sums to deduce their
limiting distribution.

\subsection{Computation of exponential sums}
For $\nu\in\bbZ/p^k\bbZ$ its ``square root'' (modulo $p^k$) is the
set
$$\Sq(\nu,p^k)=\set{x\in\bbZ/p^k\bbZ|x^2=\nu\pmod{p^k}}.$$
Note that for $\nu\neq 0\pmod {p}$ this set contains two or zero
elements, for $\nu\equiv 0\pmod {p^k}$ it contains $p^{[k/2]}$
elements (and for $\nu=p^l\tilde\nu$ with $\tilde\nu$ coprime to $p$
it contains zero or $2p^{l/2}$ elements).

\begin{prop}
For $k=2l$ even
\[E_{p^k}(\nu,\chi)=p^{l}\!\!\!\!\!\!\!\!\mathop{\sum_{x\in \Sq(\frac{2t_\chi+\nu}{\nu D},p^l)}}_{Dx^2\neq 1(p)}\!\!\!\!\!\!\!\!e_{p^k}(\nu x)\chi(\beta(x)),\]
where $t_\chi\in\bbZ/p^l\bbZ$ satisfies that
$\chi(1+p^l\sqrt{D}x)=e_{p^l}(t_\chi x)$

For $k=2l+1$ odd
\[E_{p^k}(\nu,\chi)=p^{l}\!\!\!\!\!\!\!\!\mathop{\sum_{x\in \Sq(\frac{2t_\chi+\nu}{\nu D},p^l)}}_{Dx^2\neq 1(p)}\!\!\!\!\!\!\!\!e_{p^k}(\nu x)\chi(\beta(x))\calG(x)\]
where $t_\chi\in\bbZ/p^{l+1}\bbZ$ satisfies
$\chi(1+p^l\sqrt{D}x+p^{2l}\frac{D}{2}x^2)=e_{p^{l+1}}(t_\chi x)$,
and $\calG(x)$ is the Gauss sum given by
\[\calG(x)=\sum_{y\in\bbZ/p\bbZ}e_p(f(x)y^2+g(x)y),\]
with $f(x)=\frac{2t_\chi x}{Dx^2-1}$ and $g(x)=p^{-l}(\nu
-t_\chi\frac{2}{Dx^2-1})$. (Notice that for
$x\in\Sq(\frac{2t_\chi+\nu}{\nu D},p^l)$ we have $(\nu
-t_\chi\frac{2}{Dx^2-1})\equiv 0\pmod{p^l}$, hence $p^{-l}(\nu
-t_\chi\frac{2}{Dx^2-1})$ gives a well defined residue modulo $p$).
\end{prop}
\begin{proof}
First for $k=2l$, write the sum as
\[E_{p^k}(\nu,\chi)=\sum_{x\in X(p^l)}\sum_{y\in\bbZ/p^l\bbZ}e_{p^k}(\nu (x+p^ly))\chi(\beta(x+p^ly)).\]
Replace $\beta(x+p^ly)\equiv \beta(x)(1+\frac{\beta'}{\beta}(x)p^l
y) \pmod{p^{2l}}$ to get
\[E_{p^k}(\nu,\chi)=\sum_{x\in X(p^l)}e_{p^k}(\nu x)\chi(\beta(x))\sum_{y\in\bbZ/p^l\bbZ}e_{p^l}(\nu y)\chi(1+p^l\frac{\beta'}{\beta}(x)y).\]
Differentiating $\beta(x)=\frac{\sqrt{D}x+1}{\sqrt{D}x-1}$ we get
$\frac{\beta'}{\beta}(x)=-\frac{2\sqrt{D}}{Dx^2-1}$, so that the
inner sum takes the form
\[\sum_{y\in\bbZ/p^l\bbZ}e_{p^l}(\nu y)\chi(1-p^l\sqrt{D}\frac{2y}{Dx^2-1}).\]
The map $x\mapsto 1+p^l\sqrt{D}x$ is an isomorphism of
$\bbZ/p^l\bbZ$ and $C_p(2l,l)$. Hence, for any character $\chi$ of
$C(p^{2l})$ there is $t_\chi\in\bbZ/p^l\bbZ$ such that
$\chi(1+p^l\sqrt{D}x)=e_{p^l}(t_\chi x)$. We can thus write the
inner sum as
\[\sum_{y\in\bbZ/p^l\bbZ}e_{p^l}(\nu y)e_{p^l}(-\frac{2yt_\chi}{Dx^2-1})=\sum_{y\in\bbZ/p^l\bbZ}e_{p^l}((\nu-\frac{2t_\chi}{Dx^2-1})y).\]
This sum vanishes unless $x\in \Sq(\frac{2t_\chi+\nu}{\nu D},p^l)$
in which case it is equal $p^l$.

 Now for $k=2l+1$, we start again by writing
\[E_{p^k}(\nu,\chi)=\sum_{x\in X(p^l)}\sum_{y\in\bbZ/p^l\bbZ}e_{p^k}(\nu (x+p^ly))\chi(\beta(x+p^ly)),\]
and replace
\[\beta(x+p^ly)\equiv \beta(x)(1+\frac{\beta'}{\beta}(x)p^l y +\frac{1}{2}\frac{\beta''}{\beta}(x)p^{2l}y^2) \pmod{p^{2l+1}}.\]
%  to get
% \[E_{p^k}(\nu,\chi)=\sum_{x\in X(p^l)}e_{p^k}(\frac{\nu (x)}{2})
% \chi(\beta(x))\sum_{y\in\bbZ/p^l\bbZ}e_{p^l}(\frac{\nu y}{2})\chi(1+p^l\frac{\beta'}{\beta}(x)y+p^{2l}\frac{\beta''}{2\beta}(x)y^2).\]
It easy to verify that the map $x\mapsto
1+p^l\sqrt{D}x+p^{2l}\frac{D}{2}x^2$ is an isomorphism of
$\bbZ/p^{l+1}\bbZ$ with $C_p(2l+1,l)$. Consequently, for every
character $\chi$ of $C(p^k)$, there is $t_\chi\in\bbZ/p^{l+1}\bbZ$
such that
\[\chi(1+p^l\sqrt{D}x+p^{2l}\frac{D}{2}x^2)=e_{p^{l+1}}(t_\chi x).\]

By differentiating $\beta(x)=\frac{\sqrt{D}x+1}{\sqrt{D}x-1}$
(twice), we get that

\begin{eqnarray*}
\lefteqn{\frac{\beta'}{\beta}(x)y+p^{l}\frac{\beta''}{2\beta}(x)y^2=}\\
&&=\sqrt{D}(-\frac{2(y-xp^ly^2)}{Dx^2-1})+
p^{l}\frac{D}{2}(\frac{2(y-xp^ly^2)}{Dx^2-1})^2\pmod{p^{l+1}},
\end{eqnarray*}
%&\equiv & 1-p^l\sqrt{D}\frac{2(y-xp^ly^2)}{x^2-D}+p^{2l}\frac{2D}{(x^2-D)^2}y^2\\
implying that the inner sum is of the form
 \[\sum_{y\in\bbZ/p^{l+1}\bbZ}e_{p^{l+1}}(\nu y-\frac{2t_\chi}{Dx^2-1}y+p^l\frac{2t_\chi x}{Dx^2-1}y^2).\]
This sum vanishes unless $x\in \Sq(\frac{2t_\chi+\nu}{\nu D},p^l)$.
To see this make a change of summation variable $y\mapsto y+p$ to
get that
\begin{eqnarray*}
\lefteqn{\sum_{y\in\bbZ/p^{l+1}\bbZ}e_{p^{l+1}}(\nu
y-t_\chi\frac{2}{Dx^2-1}y+p^l\frac{2t_\chi x}{Dx^2-1}y^2)=}\\
&&=e_{p^l}((\nu
-\frac{2t_\chi}{Dx^2-1}))\!\!\!\!\!\sum_{y\in\bbZ/p^{l+1}\bbZ}e_{p^{l+1}}(\nu
y-\frac{2t_\chi}{Dx^2-1}y+p^l\frac{2t_\chi x}{Dx^2-1}y^2).
\end{eqnarray*}
Now unless $\nu -\frac{2t_\chi}{Dx^2-1}\equiv 0\pmod{p^l}$ we have that $e_{p^l}((\nu -\frac{2t_\chi}{Dx^2-1}))\neq
1$, implying that the sum must vanish. For $x\in
\Sq(\frac{2t_\chi+\nu}{\nu D},p^l)$ the inner sum given by $p^l$
times the Gauss sum
\[\calG(x)=\sum_{y\in\bbZ/p\bbZ}e_p(f(x)y^2+g(x)y).\]
%with $f(x)=\frac{2t_\chi x}{Dx^2-1}$ and $g(x)=p^{-l}(\nu-t_\chi\frac{2}{Dx^2-1})$.
\end{proof}

In particular this computation implies that for most characters the
exponential sum has square root cancelation.
\begin{cor}\label{c:bound}
For any character $\chi$ with $2t_\chi\not\equiv -\nu \pmod p$ there
is $\theta=\theta(\chi,\nu)\in[0,\pi)$ such that
$E_{p^k}(\nu,\chi)=p^{k/2}\cos(\theta(\nu,\chi))$.
\end{cor}
\begin{proof}
The condition $2t_\chi\not\equiv -\nu \pmod p$ implies
$\frac{2t_\chi+\nu}{\nu D}\neq 0\pmod{p}$. Hence for any $1\leq
l\leq k$
\[\#Sq(\frac{2t_\chi+\nu}{\nu D},p^l)=\left\lbrace\begin{array}{cc}
2 &\frac{2t_\chi+\nu}{\nu D}= \square\pmod{p}\\
0 & \mbox{ otherwise}\end{array}\right.\]

Now for $k=2l$ even, recall that
\[E_{p^k}(\nu,\chi)=p^{l}\sum_{x\in\Sq(\frac{2t_\chi+\nu}{\nu D},p^l)}e_{p^k}(\nu x)\chi(\beta(x)).\]
If $\frac{2t_\chi+\nu}{\nu D}\neq \square\pmod{p}$ this sum
vanishes. Otherwise it is a sum over two elements of absolute value
$p^l=p^{k/2}$ hence indeed $E_{p^k}(\nu,\chi)=
2p^{k/2}\cos(\theta(\nu,\chi))$.

For $k=2l+1$ odd we have
\[E_{p^k}(\nu,\chi)=p^{l}\sum_{x\in \Sq(\frac{2t_\chi+\nu}{\nu D},p^l)}e_{p^k}(\nu x)\chi(\beta(x))\calG(x).\]
The condition $2t_\chi\not\equiv -\nu \pmod p$ implies that the
Gauss sum $\calG(x)$ is not a trivial sum and hence of order
$\sqrt{p}$. As before, the sum $E_{p^k}(\nu,\chi)$ either vanishes
(if $\frac{2t_\chi+\nu}{\nu D}\neq \square\pmod{p}$) or it is a sum
of two elements of absolute value $p^{l+\frac{1}{2}}=p^{k/2}$.
\end{proof}

On the other hand, if $2t_\chi\equiv -\nu \pmod{p^{2l'}}$ for some
$l'\leq \frac{k}{2}$ then the sum contains $p^{l'}$ elements and
could be much larger. Moreover, in the odd case, if $2t_\chi\equiv
-\nu\pmod{p^{l+1}}$ then $g(x)\equiv 0\pmod{p}$. Also, in this case
any $x\in\Sq(\frac{2t_\chi+\nu}{\nu D},p^l)$ satisfies $x\equiv
0\pmod{p}$ and hence also $f(x)=\frac{2t_\chi x}{Dx^2-1}\equiv
0\pmod{p}$. So that in this case the Gauss sum $|\calG(x)|=p$ rather
then $\sqrt{p}$ and the sum is even bigger. In particular we get
\begin{cor}\label{c:big}
For $\nu\in\bbZ$ and $\chi\in\hat C(p^3)$ with $t_\chi\equiv
-\nu\pmod{p^2}$,
\[|E_{p^3}(\nu,\chi)|=p^2.\]
\end{cor}

\subsection{Equidistribution of exponential sums}
We now show that as $p\to\infty$ the normalized exponential sums
$p^{-k/2}E_{p^k}(\nu,\chi)$ become equidistributed with respect to
the measure \begin{equation}\label{e:measure}
\mu(f)=\frac{1}{2}f(\frac{\pi}{2})+\frac{1}{2\pi}\int_0^{\pi}f(\theta)d\theta.
\end{equation}

For fixed $\nu$ and a character $\chi$, if $\frac{2t_\chi+\nu}{\nu
D}$ is not a square modulo $p$ then the sum $E_{p^k}(\nu,\chi)=0$
(or equivalently $\theta(\nu,\chi)=\frac{\pi}{2}$). The following
lemma shows that this happens for roughly half the characters, and
that this behavior is independent for different values of $\nu$.
\begin{lem}\label{l:square}
Fix a finite set $\bar\nu=\{\nu_1,\ldots,\nu_r\}$ of nonzero distinct integers. Then,
\[\frac{1}{p}\#\set{t\in\bbZ/p\bbZ| \forall j,\;\frac{t-\nu_j}{D\nu_j}\equiv \square\pmod{p}}=\frac{1}{2^r}+O(\frac{1}{\sqrt{p}}).\]
\end{lem}
\begin{proof}
    We can write
    \begin{eqnarray*}\lefteqn{2^r\#\set{t\in\bbZ/p\bbZ| \forall j,\;\frac{t-\nu_j}{D\nu_j}\equiv \square\pmod{p}}=}\\
   &&\quad\quad\quad\quad\quad\quad\quad\quad\quad\quad\quad\quad\quad\quad\quad=\sum_{t}\prod_{j=1}^r\left(\chi_2(\frac{t-\nu_j}{D\nu_j})+1\right),\end{eqnarray*}
    with $\chi_2$ the quadratic character modulo $p$. Now expand the
    right hand side
    \[\sum_{t}\prod_{j=1}^r\left(\chi_2(\frac{t-\nu_j}{D\nu_j})+1\right)=\sum_{J\subseteq \{1,\ldots, r\}}\sum_{t}\chi_2(\prod_{j\in
    J}\frac{t-\nu_j}{D\nu_j}).\]
    Where the sum is over all subsets $J\subseteq \{1,\ldots, r\}$.
    The contribution of the empty set is exactly $\sum_t 1=p$,
    while for nonempty $J$ we get an exponential sum of the form
    $\sum_{t}\chi_2(\prod_{j\in
    J}\frac{t-\nu_j}{D\nu_j})$. Since we assumed all $\nu_j$ are distinct,
    the polynomial $g(t)=\prod_{j\in
    J}\frac{t-\nu_j}{D\nu_j}$ is not a square and we can apply the Weil
    bounds $\sum_{t}\chi_2(g(t))=O(\sqrt{p})$ \cite{Weil48}. Consequently, we have that indeed
    \[2^r\#\set{t\in\bbZ/p\bbZ| \forall j,\;\frac{t-\nu_j}{D\nu_j}\equiv \square\pmod{p}}=p+O(\sqrt{p}).\]
\end{proof}

Next we need to show that for the rest of the characters (when the
exponential sum does not vanish) the angles  $\theta(\chi,\nu)$
become equidistributed (independently) in $[0,\pi]$. We will do that
by computing all mixed moments. However, we recall that there are
exceptional characters for which the normalized exponential sums are
not bounded causing the moments to blow up. For that reason we first
restrict ourself to a set of ``good" characters (of limiting density
one) for which the sums are bounded and only then we calculate the
moments.

Fix a finite set of $r$ nonzero distinct integers
$\bar\nu=\{\nu_1,\ldots,\nu_r\}$, and define the set of ``good"
characters to be
\[S_{p^k}(\bar\nu)=\set{\chi\in\hat C(p^k)|\forall j, 2t_\chi \not\equiv -\nu_j\pmod{p}},\]
where $t_\chi$ is determined by $\chi$ as above. Then for any
character $\chi\in S_{p^k}(\bar\nu)$, we can write
$E_{p^k}(\nu_j,\chi)=p^{k/2}\cos(\theta(\nu_j,\chi))$ with
$\theta(\nu_j,\chi)\in[0,\pi)$. Furthermore, for any  $\nu_j$ there
are precisely $p^{k-2}(p\pm 1)$ characters with
$2t_\chi\equiv\nu_j\pmod{p}$ (this is the size of the kernel of the
restriction map from $\hat C(p^k)$ to $\hat C_p(k,1)$). Hence,
$\frac{|S_{p^k}(\bar\nu)|}{p^k}=1+O(\frac{1}{p})$ and the set
$S_{p^k}(\bar\nu)$ is of (limiting) density one inside $\hat
C(p^k)$.

Before we proceed to calculate the moments we will need to set some
notations. For any $k$ define the set
\[Y(p^k,\bar \nu)=\set{\bar{x}\in (X(p^k))^r|\nu_1(Dx_1^2-1)=\nu_j(Dx_j^2-1),\;\forall 2\leq j\leq r}\]
For every fixed set of integers $\bar n=\set{n_1,\ldots,n_r}$ let
\[Y_0(p^k,\bar\nu,\bar n)=\set{\bar{x}\in Y(p^k)|\prod_j \beta(x_j)^{n_j}\equiv 1\pmod{p^k}}.\]
For notational convenience we will sometimes use the notation
$S_{p^k}$, 
$Y(p^k)$, $Y_0(p^k)$ where the dependence on $\bar\nu$ and
$\bar n$ is implicit. We will also denote by $Y'(p^k)$ (respectively
$Y'_0(p^k)$) the elements of $Y(p^k)$ (respectively $Y_0(p^k)$) with
all $x_j\neq 0\pmod{p}$.

\begin{lem}\label{l:count1}
As $p\to\infty$, the number of points in $Y'(p^k)$ satisfy
\[\# Y'(p^k)=p^k+O(p^{k-\frac{1}{2}})\]
\end{lem}
\begin{proof}
For any $t\in\bbZ/p^k\bbZ$ satisfying $\forall
j,\;t\neq\nu_j\pmod{p}$ we have that
\[\#\set{\bar{x}|\forall j,\;\nu_j(Dx_j^2-1)=t}=\left\lbrace\begin{array}{cc}
2^r & \forall j,\;\frac{t-\nu_j}{D\nu_j}\equiv \square\pmod{p}\\
 0 & \mbox{ otherwise}
\end{array}\right.
\]
On the other hand if $t\equiv\nu_j\pmod{p}$ for some $j$, then
\[\#\set{\bar{x}|\forall j,\;\nu_j(Dx_j^2-1)=t}\leq 2^{r}p^{k-1}\]
(as there are at most two possibilities for $x_i$ with $i\neq j$ and
at most $2p^{k-1}$ possibilities for $x_j$). We thus have
\begin{eqnarray*}\lefteqn{\# Y(p^k) = \sum_{t\in\
 (\bbZ/p^k\bbZ)^*}\#\set{\bar{x}\in(\bbZ/p^k\bbZ)^r|\nu_j(Dx_j^2-1)=t}}\\
&&= 2^rp^{k-1}\#\set{t\in(\bbZ/p\bbZ)^*| \forall
j,\;\frac{t-\nu_j}{D\nu_j}\equiv \square\pmod{p}}+O(p^{k-1}).
\end{eqnarray*}
Also note that  $\# Y'(p^k)=\# Y(p^k)+O(p^{k-1})$. To conclude the
proof we use the estimate
\[2^r\#\set{t\in\bbZ/p\bbZ| \forall j,\;\frac{t-\nu_j}{D\nu_j}\equiv \square\pmod{p}}=p+O(\sqrt{p}),\]
from lemma \ref{l:square}.

\end{proof}

\begin{lem}\label{l:count2}
As $p\to\infty$, the number of points in $Y'_0(p^k)$ satisfy
$$\# Y'_0(p^k)=O(p^{k-1}).$$
\end{lem}
\begin{proof}
To prove this bound we will show that there is a nonzero polynomial
$F(t)$ with integer coefficients such that for any $\bar{x}\in
Y_0(p^k)$, with $b=\frac{\nu_1(Dx_1^2-1)}{2}$ we have
$F(b^{-1})\equiv 0\pmod{p^k}$ (recall that for $x_1\in X(p^k)$ we
have $Dx_1^2\neq 1\pmod{p}$ and hence $b\neq 0\pmod{p}$ is
invertible). This would imply that $b$ can take at most $\deg F$
values modulo $p$, implying that $\# Y_0(p^k)\leq
2^r\deg(F)p^{k-1}$.

Now to define $F$, consider the formal polynomial in the variables
$\beta_1^{\pm 1},\ldots\beta_r^{\pm 1}$ given by
\[G(\beta_1,\ldots\beta_r)=\prod_{\sigma\in\{\pm 1\}^r}\left(\prod_{j=1}^r\beta_j^{\sigma_jn_j}-1\right).\]
Recall that if a polynomial in two variables $x,y$ is symmetric
under permutation then it can be written as a polynomial in the
symmetric polynomials $\sigma_1=x+y,\sigma_2=xy$ (see e.g.,
\cite[Chapter 6]{Bourbaki90}). The polynomial $G$ is symmetric under
any substitution $\beta_j\mapsto\beta_j^{-1}$ and hence there is
another polynomial $\tilde{F}$ in $r$ variables with integer
coefficients, satisfying
 \[G(\beta_1,\ldots,\beta_r)=\tilde{F}(\beta_1+\beta_1^{-1},\ldots,\beta_r+\beta_r^{-1}).\]

Define the polynomial $F(t)=\tilde{F}(2+\nu_1t,\ldots,2+\nu_r t)$.
For any $x_1,\ldots,x_r$ with $x_j^2=\frac{2b-\nu_j}{\nu_j
D}\pmod{p^k}$ we have
$\beta(x_j)+\beta(x_j)^{-1}=2+\nu_jb^{-1}\pmod{p^k}$ (recall
$\beta(x)=\frac{\sqrt{D}x+1}{\sqrt{D}x-1}$). Hence,
\[G(\beta(x_1),\ldots,\beta(x_r))=\tilde{F}(2+\nu_1b^{-1},\ldots,2+\nu_rb^{-1})=F(b^{-1}).\]
Now, if in addition $\beta(x_1)^{n_1}\cdots\beta(x_r)^{n_r}\equiv
1\pmod {p^k\frakO}$ then  indeed
$F(b^{-1})=G(\beta(x_1),\ldots,\beta(x_r))\equiv 0\pmod{p^k}$.

It remains to show that $F(t)$ is not the zero polynomial. To do
this, we think of it as a complex valued polynomial, and note that
for it to be identically zero there has to be some choice of signs
$\sigma\in\{\pm 1\}^r$ so that the function
\[G_\sigma(t)=\prod_{j=1}^r \beta(\sqrt{\frac{2t+\nu_j}{\nu_j D}})^{\sigma_j n_j}\]
satisfies $G_\sigma(t)\equiv1$. Assume that there is such a choice
$\sigma$, so the derivative $G_\sigma'(t)$ must also vanish. But we
have
\[G'(t)=-G_\sigma(t)\sum_{j=1}^r\sigma_jn_j\sqrt{\frac{\nu_j}{t^2(2t+\nu_j)}},\]
so as $t\to-\frac{\nu_1}{2}$ the term
$\sqrt{\frac{\nu_1}{t^2(2t+\nu_1)}}$ blows up while the rest of the
terms remain bounded (recall that all $\nu_j$ are different). In
particular $G'_\sigma(t)$ is not identically zero.
\end{proof}
\begin{rem}
The bound $\# Y'_0(p^k)=O(p^{k-1})$ is probably not optimal. Notice
that if the polynomial $F(t)$ defined above is separable (i.e., if
it has no multiple roots) then there are at most $\deg F$ solutions
to $F(t)\equiv 0\pmod{p^k}$ and the corresponding bound would be $\#
Y'_0(p^k)=O(1)$.
\end{rem}

We now preform the moment calculation establishing the limiting
distribution of the exponential sums (when running over characters
in $S_{p^k}$).
\begin{prop}\label{p:equi}
Let $\mu$ be as in (\ref{e:measure}) and let $g\in C([-1,1]^r)$ be
any continuous function then
\begin{eqnarray*}
\lefteqn{\lim_{p\to\infty}\frac{1}{p^k}\sum_{\chi\in S_{p^k}}g(\cos(\theta(\nu_1,\chi),\ldots,\cos(\theta(\nu_r,\chi))=}\\
&&\int_{[0,\pi]^d}g(\cos(\theta_1),\ldots,\cos(\theta_r))d\mu(\theta_1)\cdots
d\mu(\theta_r).
\end{eqnarray*}
\end{prop}
\begin{proof}
We will give the proof for $k=2l$ even, the odd case is analogous.
Since we can always approximate the function $g$ by polynomials, it
is sufficient to show this holds for all monomials of the form
$$g(x)=(2x_1)^{m_1}\cdots (2x_r)^{m_r}.$$
We thus need to show that
\begin{eqnarray*}
\lim_{p\to\infty}\frac{1}{p^k}\sum_{\chi\in
S_{p^k}}\prod_j(2\cos(\theta(\nu_j,\chi)))^{m_j}=\prod_{j}\int_{[0,\pi]}(2\cos(\theta))^{m_j}d\mu(\theta).
\end{eqnarray*}
With out loss of generality we can also assume that all the $m_j$
are nonzero (since $\mu$ is a probability measure, if $m_j=0$ then
the corresponding factor is $1$ and we can consider the same problem
for $r-1$ instead of $r$). In this case the right hand side is given
by
\begin{eqnarray*}
 \prod_j\left(\int_0^\pi \left(2\cos(\theta)\right)^{m_j}d\mu(\theta)\right)=\prod_j\left(\frac{1}{2}\int_0^\pi \left(2\cos(\theta)\right)^{m_j}\frac{d\theta}{\pi}\right).
\end{eqnarray*}
The integral in each factor is
$\frac{1}{2}\begin{pmatrix}m_j\\n_j\end{pmatrix}$ for $m_j=2n_j$
even and it is zero otherwise.

Now fix a character $\chi\in S_{p^k}$ and let
$t_\chi\in\bbZ/p^l\bbZ$ as above. If $\Sq(\frac{2t_\chi+\nu_j}{\nu_j
D},p^l)=\emptyset$ then
$2\cos\theta(\nu_j,\chi)=p^{-k/2}E_{p^k}(\nu_j,\chi)=0$. Otherwise,
\[2\cos(\theta(\nu_j,\chi))=p^{-{k/2}}E_{p^k}(\nu_j,\chi)=2\Re(e_{p^k}(\frac{\nu_j x_j}{2})\chi(\beta(x_j))),\]
with $x_j\in\Sq(\frac{2t_\chi+\nu_j}{\nu_j D},p^l)$ (recall that for
$\chi\in S_{p^k}$ we know $2t_\chi+\nu\neq 0\pmod{p}$). Hence, the
only contributions to the sum
\[\sum_{\chi\in S_{p^k}}\prod_{j} \left(2\cos\theta(\nu_j,\chi)\right)^{m_j},\]
comes from characters $\chi$ such that for all $j$ there is
$x_j\in\Sq(\frac{2t_\chi+\nu_j}{\nu D},p^l)$ (equivalently, there is
$x_j\in(\bbZ/p^l\bbZ)^*$ satisfying $\nu_j(Dx_j^2-1)\equiv
2t_\chi\pmod{p^l}$). Also note that if we multiply $\chi$ by any
character that is trivial on $C_p(k,l)$ this does not change
$t_\chi$. Let $\hat C^{(l)}(p^k)$ be the group of characters that
are trivial on $C_p(k,l)$, and for any $b\in\bbZ/p^l\bbZ$ let
$\chi_b\in \hat C(p^k)$ be a representative of $\hat C(p^k)/\hat
C^{(l)}(p^k)$ with $t_{\chi_b}=b$. We thus have that
\begin{eqnarray*}\lefteqn{\frac{1}{p^k}\sum_{\chi\in S_{p^k}}\prod_{j} \left(2\cos\theta(\nu_j,\chi)\right)^{m_j}=}\\
&&=\frac{1}{2^rp^k}\sum_{\bar{x}\in Y'(p^l)}\sum_{\chi\in \hat
C^{(l)}(p^k)}
\prod_{j}\left(2\cos\theta(\nu_j,\chi\chi_b)\right)^{m_j}
\end{eqnarray*}
where $b=b(\bar{x})=\frac{\nu_1(Dx_1^2-1)}{2}$.

Now use the formula,
\[(2\cos(\theta))^m=\sum_{n=0}^m\begin{pmatrix} m\\n\end{pmatrix}\cos((m-2n)\theta).\]
The main contribution comes from the terms where in each factor
$m_j-2n_j=0$. This vanishes unless all $m_j$ are even in which case
it is given by
\[\frac{1}{2^rp^k}\sum_{\bar{x}\in Y'(p^l)}\sum_{\chi\in \hat C^{(l)}(p^k)}
\prod_{j}\begin{pmatrix} m_j\\
n_j\end{pmatrix}=\prod_{j}\frac{1}{2}\begin{pmatrix} m_j\\
n_j\end{pmatrix}+O(\frac{1}{\sqrt{p}}).\] where we used Lemma
\ref{l:count1} to get that $\# Y'(p^l)\cdot\# \hat
C^{(l)}(p^k)=p^k+O(p^{k-\frac{1}{2}})$.

It thus remains to bound the rest of the terms, which is reduced to
the vanishing (in the limit $p\to\infty$) of the sums
\[\frac{1}{p^k}\sum_{\bar{x}\in Y'(p^l)}\sum_{\chi\in \hat C^{(l)}(p^k)}
\prod_{j}\cos(n_j\theta(\nu_j,\chi\chi_b)),\] for any nonzero
integers $\{n_1,\ldots, n_r\}$.

For any $\bar{x}\in Y'(p^l)$ we have that
\[\cos(n_j\theta(\nu_j,\chi\chi_b)=2\Re(e_{p^k}(\frac{n_j \nu x_j}{2})\chi\chi_b(\beta(x_j)^{n_j}))\]
with $b=\frac{\nu_1(Dx_1^2-1)}{2}$. When expanding the product
$\prod_{j}\cos(n_j\theta(\nu_j,\chi\chi_b))$ we get a sum over $2^r$
terms, each of the form
$$e_{p^k}(\sum_{j=1}^r
\frac{\pm n_j \nu_j x_j}{2})\chi\chi_b(\prod_{j=1}^r\beta(x_j)^{\pm
n_j}).$$ We thus need to bound the exponential sum coming from each
term. We will now bound the corresponding sum
\[\frac{1}{p^k}\sum_{\bar{x}\in Y'(p^l)}\sum_{\chi\in \hat
C^{(l)}(p^k)}e_{p^k}(\sum_{j=1}^r\frac{ n_j \nu_j
x_j}{2})\chi\chi_b(\prod_{j=1}^r\beta(x_j)^{ n_j}).\] (the same
bound obviously holds when changing any $n_j$ to $-n_j$). Now, rewrite
this sum as
\[\frac{1}{p^k}\sum_{\bar{x}\in Y'(p^l)}e_{p^k}
(\sum_{j=1}^r\frac{n_j \nu_j
x_j}{2})\chi_b(\prod_{j=1}^r\beta(x_j)^{n_j})\sum_{\chi\in \hat
C^{(l)}(p^k)}\chi(\prod_{j=1}^r\beta(x_j)^{n_j}),\] and note that
the inner sum vanishes unless $\prod_{j=1}^r\beta(x_j)^{n_j}\equiv
1\pmod {p^l}$ in which case it is equal $\#\hat
C^{(l)}(p^k)=p^{k-l}$. We can thus rewrite this sum as
\[\frac{1}{p^l}\sum_{\bar{x}\in Y_0'(p^l)}e_{p^k}
(\sum_{j=1}^r\frac{n_j \nu_j
x_j}{2})\chi_b(\prod_{j=1}^r\beta(x_j)^{n_j}).\] which is trivially
bounded by $p^{-l}\# Y'_0(p^l)=O(\frac{1}{p})$ (Lemma
\ref{l:count2}).
\end{proof}
\begin{rem}
The above proof also gives the rate at which the fluctuations of the
normalized exponential sums approach their limiting distribution. If
one takes the test function $g$ in Proposition \ref{p:equi} to be
smooth then the rate of convergence is $O(\frac{1}{\sqrt{p}})$. This
rate comes from the bound on the error term in Lemma \ref{l:count1}
which seems to be a sharp bound.
\end{rem}
\section{Back to Matrix Elements} \label{s:equi}
We can now deduce Theorems \ref{t:slowdecay} and \ref{t:main} from
Theorem \ref{t:form1} and the analysis of the exponential sums.

\begin{proof}[\bf{Proof of Theorem \ref{t:slowdecay}}]
Let $f(x)=e^{2\pi i n\cdot x}$ be any elementary observable. Take
$N=p^3$ to be a prime cubed. Then by Corollary \ref{c:big} there is
a character satisfying $|E_{p^3}(\frac{Q(n)}{2},\chi)|=p^2$. Let
$\psi$ be a Hecke eigenfunction corresponding to $\chi$, then by
Theorem \ref{t:form1} we get
\[|\langle \Op_{N}(f)\psi,\psi\rangle|=\frac{1}{\# C(p^3)}E_{p^3}(\frac{Q(n)}{2},\chi)=\frac{1}{p\pm 1}\gg N^{-1/3}.\]
\end{proof}

\begin{proof}[\bf{Proof of Theorem \ref{t:main}}]
Let $f$ be a trigonometric polynomial and write
\[f=\sum_{|n|\leq R}\hat{f}(n)e(n\cdot x),\]
for some fixed $R>0$. Let
$\{\nu_1,\ldots,\nu_r\}=\set{Q(n)|0<|n|\leq R}$, and consider the
random variable
\[Y_f=2\sum_{j=1}^r f^\#(\nu_j)\cos(\theta_j).\]
with $\theta_j$ chosen independently from $[0,\pi)$ with respect to
$\mu$. We need to show that as $p\to\infty$ the limiting
distribution of $F_j^{(p^k)}$ is that of $Y_f$.

For any character $\chi$ of $C(p^k)$ consider the weighted sum of
the corresponding exponential sums
\[F^{(p^k)}_\chi=\sum_{j=1}^r  f^\#(\nu_j)p^{-k/2}E_{p^k}(\frac{\nu_j}{2},\chi\chi_0).\]
By Proposition \ref{p:equi}, as $p\to\infty$ the limiting
distribution of $F^{(p^k)}_\chi$ as $\chi$ runs through $S_{p^k}$
(hence, also as $\chi$ runs through the whole group of characters)
is that of $Y_f$.

Now, for $p$ sufficiently large (i.e., $p>\max\{\nu_j\}$) and
$\chi_j\in \hat{C}_0(p^k)$, we have
\[F_j^{(p^k)}=\sum_{j=1}^r  f^\#(\nu_j)\frac{p^{k/2}}{\# C(p^k)}E_{p^k}(\frac{\nu_j}{2},\chi_j\chi_0).\]
If we further assume that $\chi_j\chi_0\in S_{p^k}(\bar\nu)$ then
$|E_{p^k}(\frac{\nu_j}{2},\chi\chi_0)|\leq2p^{k/2}$, and hence
%$\frac{p^{k/2}}{\# C(p^k)}=p^{-k/2}(1+O(\frac{1}{p}))$ we can write
\[F_j^{(p^k)}=F^{(p^k)}_{\chi_j}+O(\frac{1}{p}).\]
The set of characters $\{\chi_j\in \hat{C}_0(p^k)|\chi_j\chi_0\in
S_{p^k}\}$ is again of density one, hence, the limiting distribution
of $F_j^{(p^k)}$ is the same as of $F^{(p^k)}_{\chi_j}$ concluding
the proof.
\end{proof}

%%-------------------------
%%GATHER{/home/member/kelmerdu/My_Documents/Tex/Bib/Mybib}   % For Gather Purpose Only
%%\bibliographystyle{amsplain}
%%\bibliography{/home/member/kelmerdu/My_Documents/Tex/Bib/Mybib.bib}
%%*********************************************************************
\def\cprime{$'$}
\providecommand{\bysame}{\leavevmode\hbox
to3em{\hrulefill}\thinspace}
\providecommand{\MR}{\relax\ifhmode\unskip\space\fi MR }
% \MRhref is called by the amsart/book/proc definition of \MR.
\providecommand{\MRhref}[2]{%
  \href{http://www.ams.org/mathscinet-getitem?mr=#1}{#2}
} \providecommand{\href}[2]{#2}

%********************************************************************
\end{document}